\documentclass[12pt]{article}
\usepackage{amssymb,amsthm,vmargin,url}
\usepackage{amsmath} 
\usepackage{graphicx}

\usepackage{pxfonts}
\usepackage[footnotesize]{caption}

\setpapersize{USletter}
\setmargnohfrb{1.0in}{1.0in}{1.0in}{1.0in}

\AtBeginDocument{
\headsep=30pt
}

\newtheorem{Th}{Theorem}
\newtheorem{Lemma}[Th]{Lemma}
\newtheorem{Def}[Th]{Definition}

\newcommand{\Z}{\mathbb Z}
\newcommand{\N}{\mathbb N}
\newcommand{\R}{\mathbb R}
\newcommand{\K}{{\cal K}}
\newcommand{\M}{{\cal F}}
\renewcommand{\L}{\mathfrak L}

\newcommand{\0}{{\mathbf 0}}

\newcommand{\eqdef}{\mathbin{\mathop=\limits^{\rm def}}}

\let\oldsum\sum
\renewcommand{\sum}{\oldsum\limits}
\let\oldint\int
\renewcommand{\int}{\oldint\limits}
\renewcommand{\l}{\ell}

\newcommand{\ceil}[1]{\left\lceil #1\right\rceil}
\newcommand{\floor}[1]{\left\lfloor #1\right\rfloor}
\newcommand{\ls}{\left(}\newcommand{\rs}{\right)}

\newcommand{\bea}{\begin{array}}
\newcommand{\eea}{\end{array}}
\newcommand{\beq}{\begin{equation}}
\newcommand{\eeq}{\end{equation}}

\newcommand{\card}[1]{\mathop{\#}#1}

\title{On the number of two-dimensional threshold functions}
\author{Max A. Alekseyev\\
$\langle\mbox{\tt maxal@cse.sc.edu}\rangle$\\
\\
Department of Computer Science and Engineering\\
University of South Carolina, U.S.A.\\
}

\begin{document}
\maketitle\thispagestyle{empty}

A two-dimensional threshold function of $k$-valued logic can be viewed as coloring of the points 
of a $k\times k$ square lattice into two colors such that there exists a straight line separating points
of different colors. For the number of such functions only asymptotic bounds are known. 
We give an exact formula for the number of two-dimensional threshold functions and derive more accurate asymptotics.

\newpage

\section{Introduction}

A function $f$ of $n$ variables of $k$-valued logic is called a \emph{threshold function} if it takes 
two values $0$ and $1$ and there exists a hyperplane separating the pre-images $f^{-1}(0)$ and $f^{-1}(1)$. 
Threshold functions have been studied from the perspectives of electrical engineering \cite{Muroga71}, neural networks \cite{Hassoun95}, 
combinatorial geometry \cite{Abelson77,Ojha00,Irmatov01}, and learning theory \cite{ShevZol98,ZolShev99}.

Computing the number $P(k,n)$ of $n$-dimensional $k$-valued threshold functions turns out to be a hard problem, even in the case of $k=2$.
The number $P(2,n)$, corresponding to $n$-dimensional boolean threshold functions, was studied in a number of 
publications \cite{Winder65,Muroga70,Ojha00,Zunic04}. 
In spite of many efforts, the exact values of $P(2,n)$ are known only for $n\leq 8$ (sequence \emph{A000609} in \cite{OEIS}).
The asymptotics of $P(2,n)$ was found in \cite{Zuev89,Zuev91}. 
Computing $P(k,n)$ for $k>2$ appears to be even a harder problem. 
Known results on the number $P(k,n)$ for $k>2$ are mostly of an asymptotic nature~\cite{Irmatov98,Irmatov01}. 
A particular case of two-dimensional threshold functions (i.e., for $n=2$) was studied in~\cite{ShevZol98}. 
Shevchenko~\cite{Shev} states the following asymptotic bounds
$$
\frac{3}{8\pi^2}k^4 \lesssim P(k,2) \lesssim \frac{6}{\pi^2}k^4.
$$
In this paper we prove (in Theorem~\ref{Total}) an exact formula for $P(k,2)$:\footnote{The values of $P(k,2)$ and $\frac{1}{2}P(k,2)$
for $k=1,2,3,\dots$ form respectively the sequences \emph{A114146} and \emph{A114043} in \cite{OEIS}.}
\begin{equation}\label{FormulaPk2}
P(k,2) = (2k-1)^2 + 1 + 4 \mathop{\sum_{i=1}^{k-1}\sum_{j=1}^{k-1}}\limits_{\gcd(i,j)=1} (k-i)(k-j)
\end{equation}
and derive (in Theorem~\ref{asympN}) a more accurate asymptotics:
\begin{equation}\label{AsympPk2}
P(k,2) = \frac{6}{\pi^2} k^4 + O(k^3\ln k).
\end{equation}

The paper is organized as follows. In Section~\ref{SecDef} we give a rigorous
definition and examples of two-dimesional threshold functions. 
The exact formula for $P(k,2)$ is obtained in Section~\ref{Sexact}. 
The asymptotics for $P(k,2)$ is derived in Section~\ref{Sasympt}. 
Finally, in Section~\ref{SecTeaching} we discuss the connection
to teaching sets and pose a related open problem.


\section{Two-dimensional threshold functions}\label{SecDef}

We consider the problem of finding $P(k,2)$ and its asymptotics in a slightly more general form, allowing 
the arguments of two-dimensional threshold functions take different numbers of integer values.
The precise definitions follow:

\begin{Def}
Let $m,\ n\in\N$ be positive integers and $\K\eqdef [0,m]\times[0,n]$ be a rectangle on Euclidean plane
bounded by the lines $x=0$, $x=m$, $y=0$, and $y=n$. The four sides of $\K$ form its \emph{border} and 
a point in $\K$ is called \emph{internal} if it does not belong to the border of $\K$.
We further define $\K_0\eqdef\K\cap\Z^2$ as the set of all integer points in $\K$. 
\end{Def}

\begin{Def} A (two-dimensional) \emph{threshold function} on $\K_0$ is a function 
$f:\K_0\to\{0,1\}$ such that there exists a line $\l = ax+by+c \equiv 0$ satisfying
$$
f(x,y)=0\qquad\Longleftrightarrow\qquad ax+by+c\leq 0.
$$
We say that the line $\l$ \emph{defines} the threshold function $f$.
\end{Def}

\begin{figure}[t]
\begin{center}
\includegraphics[scale=0.6]{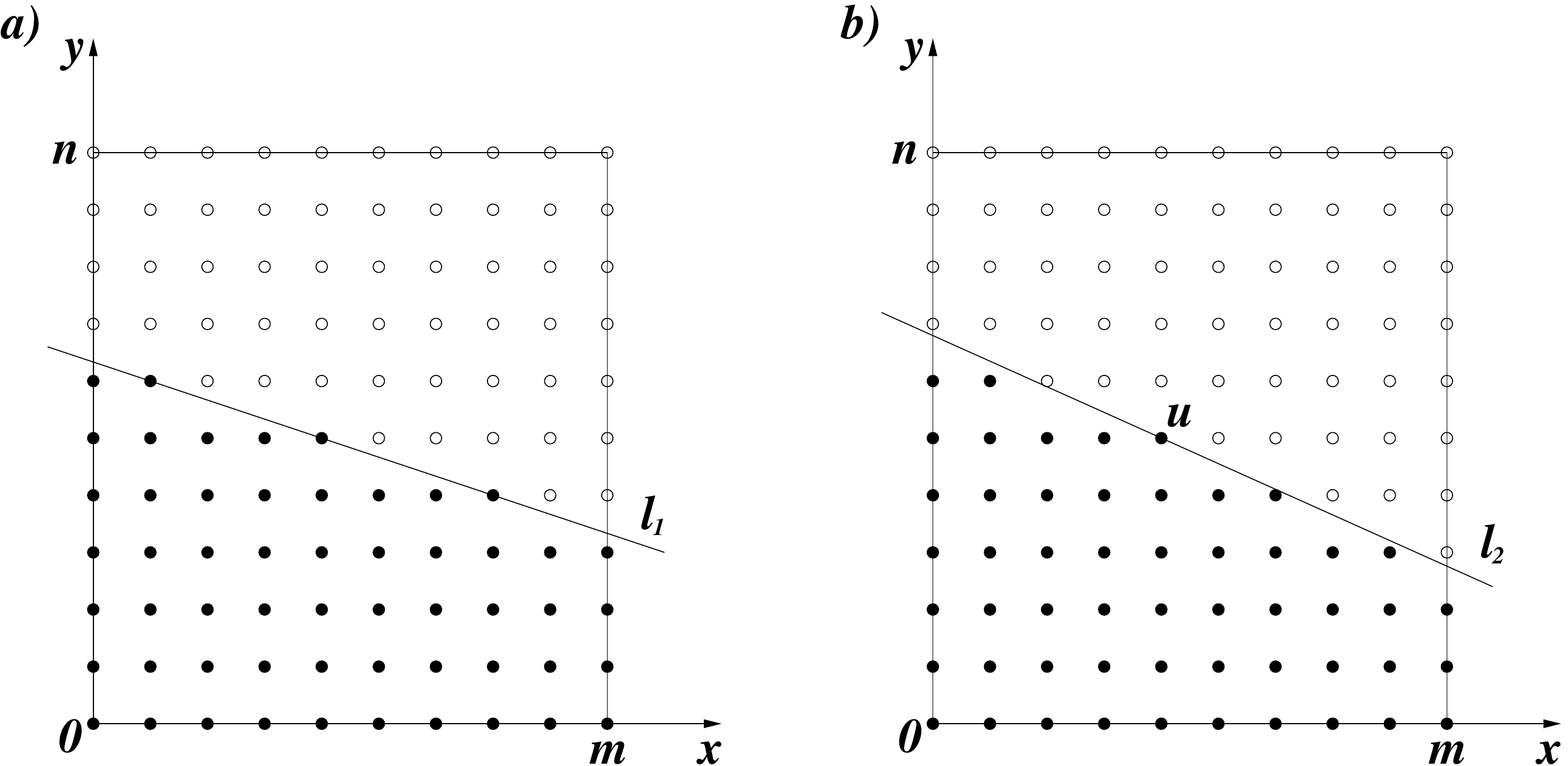}
\caption{Examples of threshold functions defined by the lines $\l_1$ and $\l_2$ (the filled dots correspond to zero values).
a) The line $\l_1$ defines a stable threshold function.
b) The line $\l_2$ defines an unstable threshold function with the vertex at $u$.}
\label{Fig0}
\end{center}
\end{figure}

Examples of threshold functions are given in Fig.~\ref{Fig0}.

Let $N(m,n)$ be the number of all threshold functions on $\K_0$.
Our goal is to find an exact formula and an asymptotics for $N(m,n)$.
That will immediately imply similar results for $P(k,2)$ since
\begin{equation}\label{FormulaPN}
P(k,2)=N(k-1,k-1).
\end{equation}


\section{Exact formula for $N(m,n)$}\label{Sexact}

In this section we prove the following theorem, which together with \eqref{FormulaPN} implies formula \eqref{FormulaPk2}.

\begin{Th}\label{Total} The total number of two-dimensional threshold functions is
$$
N(m,n) = (2m+1)(2n+1) + 1 + 4V(m,n),
$$
where 
$$
V(m,n) \eqdef \mathop{\sum_{i=1}^{\ceil{m}}\sum_{j=1}^{\ceil{n}}}\limits_{\gcd(i,j)=1} (m+1-i)(n+1-j).
$$
\end{Th}

\subsection{Preliminary results}

All threshold functions fall into two classes, depending on the value of $f(0,0)$.
Between these two classes there is a natural bijection $f\mapsto 1-f$.
Let $\M$ be a class of threshold functions with $f(0,0)=0$, excluding the zero function.
Then the total number of threshold functions equals $2(|\M|+1)$.

With each line $\l = ax+by+c \equiv 0$ we associate a set of zeros of 
a threshold function defined by $\l$:
$$
M(\l) = M(a,b,c) = \{ (x,y)\in\K_0\mid ax+by+c\leq 0\}.
$$
If the line $\l$ defines a function from $\M$, then $c\leq 0$. For the rest, we assume that
this inequality always holds.

We note that $M(\l)$ is well-defined only for a line $\l = ax+by+c \equiv 0$ with $c\ne 0$ (a \emph{regular} line); however,
a line with $c=0$ (a \emph{singular} line) generally defines two threshold functions corresponding to $M(a,b,0)$ and $M(-a,-b,0)$.
For a singular line $\l'$ resulting from continuous motion of a regular line $\l$, 
we assume that $\l'$ defines a function equal the limit of the function defined by $\l$.

\begin{Def} Call a line $\l = ax+by+c \equiv 0$ \emph{horizontal}, if $a=0$;
\emph{vertical}, if $b=0$; and \emph{inclined}, if $a\ne0$ and $b\ne0$.
An inclined line is called \emph{positive} or \emph{negative} depending on
the sign of its slope $\frac{-a}{b}$.
\end{Def}

\begin{Def} Lines $\l_1$ and $\l_2$ are \emph{equivalent} $(\l_1\sim \l_2)$, if
$M(\l_1)=M(\l_2)$.
In other words, two lines are equivalent if they define 
the same threshold function.
\end{Def}

We prove (non-)equivalence of the lines $\l_1$ and $\l_2$ using
(non-)emptiness of the symmetric difference $M(\l_1)\bigtriangleup M(\l_2)$, 
or both set differences $M(\l_1)\setminus M(\l_2)$ and $M(\l_2)\setminus M(\l_1)$.

\begin{Lemma}\label{Lprev} For any line $\l$ defining a threshold function from $\M$,
there exists an equivalent line $\l'$ passing through at least one point from $\K_0$.
\end{Lemma}

\begin{proof} Let the line $\l=ax+by+c\equiv 0$ define some threshold function from
$\M$. If $\l$ does not contain points from $\K_0$, then
increasing $c$ (i.e., translating $\l$ towards the origin) we will find a line
$\l'=ax+by+c'\equiv 0$ such that $\l'$ passes through at least one point from
$\K_0$, and there are no points from $\K_0$ between $\l'$ and $\l$.
Then $\l\sim \l'$.
\end{proof}

Denote by $\L$ the set of all lines that define functions from $\M$ and pass
through at least one point from $\K_0$. Lemma~\ref{Lprev} implies that
every function from $\M$ is defined by some line from $\L$. 

\begin{Lemma}\label{l2q} Let $\l_1,\ \l_2\in\L$. If $\l_1\sim \l_2$ and $\l_1\ne \l_2$,
then there exists a point $q\in \l_1\cap \l_2\cap\K$.
\end{Lemma}
\begin{proof}
Assume that the lines $\l_1$ and $\l_2$ do not have a common point within the rectangle $\K$.
Since $\l_1,\ \l_2\in\L$, there exist points $u\in \l_1\cap\K_0$ and $v\in \l_2\cap\K_0$.
Then either $u\in M(\l_1)\setminus M(\l_2)$ or $v\in M(\l_2)\setminus M(\l_1)$,
a contradiction to $\l_1\sim \l_2$.
\end{proof}

\begin{Def}
A line $\l\in\L$ is called \emph{stable} if it passes through at least two points from $\K_0$.
\end{Def}

An example of a stable line is given in Fig.~\ref{Fig0}a.

\begin{Lemma}\label{l12l0} Let $\l_1\sim \l_2$ be equivalent lines passing
through points $u\in\K_0$ and $v\in\K_0$ respectively. If $u\ne v$
then there exists a stable line $\l_0$ such that $u\in \l_0$ and $\l_1\sim \l_0\sim \l_2$.
\end{Lemma}
\begin{proof} If $\l_1$ is stable, then the statement is 
trivial for $\l_0=\l_1$. Assume that $\l_1$ is not stable. 

Lemma \ref{l2q} implies $\l_1\cap \l_2=q\in\K$.
Since $\l_1$ is not stable, $q\ne v$. 
On the other hand, if $q=u$ then $\l_0=\l_2$ proves the lemma. Hence, assume $q\ne u$.

Consider a family of lines $\l(t)$ passing through the points $u$ and $q+(v-q)t$
for $t\in[0,1]$. Note that $\l(0)=\l_1$.
Let $t_0>0$ be a minimal value of $t$ such that line $\l(t)$ passes through
a point from $\K_0$ different from $u$. Define $\l_0=\l(t_0)$.

We will show that the set $M(\l(t))$ does not change as the parameter $t$ changes from $0$ to $t_0$. 
By the construction, $M(\l')=M(\l_1)$ for any intermediate line $\l'=\l(t)$, $t\in(0,t_0)$ and
all the points from $\l_1\cap\K_0$, except $u$, lie on the same side of $\l'$ as the point $v\in M(\l_1)=M(\l_2)$. 
Hence, $M(\l')=M(\l_1)$.
On the other hand, for any point $w\in \l_0\cap\K_0$, it is true that $w\in M(\l_1)=M(\l_2)$. 
Indeed, if $w\not\in M(\l_1)=M(\l_2)$, then 
the point $w$ must lie on a ray of $\l_0$ starting at $u$ that crosses the line $\l_2$.
But then the points $w$ and $v$ lie at the same side of line $\l_1$, 
and thus $w\in M(\l_1)$ which contradicts the assumption $w\not\in M(\l_1)=M(\l_2)$.

Therefore, $M(\l_0)=M(\l(t_0))=M(\l(0))=M(\l_1)$, i. e., $\l_0\sim \l_1$.
\end{proof}

\begin{Lemma}\label{uK0} If equivalent stable lines $\l_1$ and $\l_2$ intersect at a point from $\K_0$, then $\l_1 = \l_2$.
\end{Lemma}

\begin{figure}[t]
\begin{center}
\includegraphics[scale=0.6]{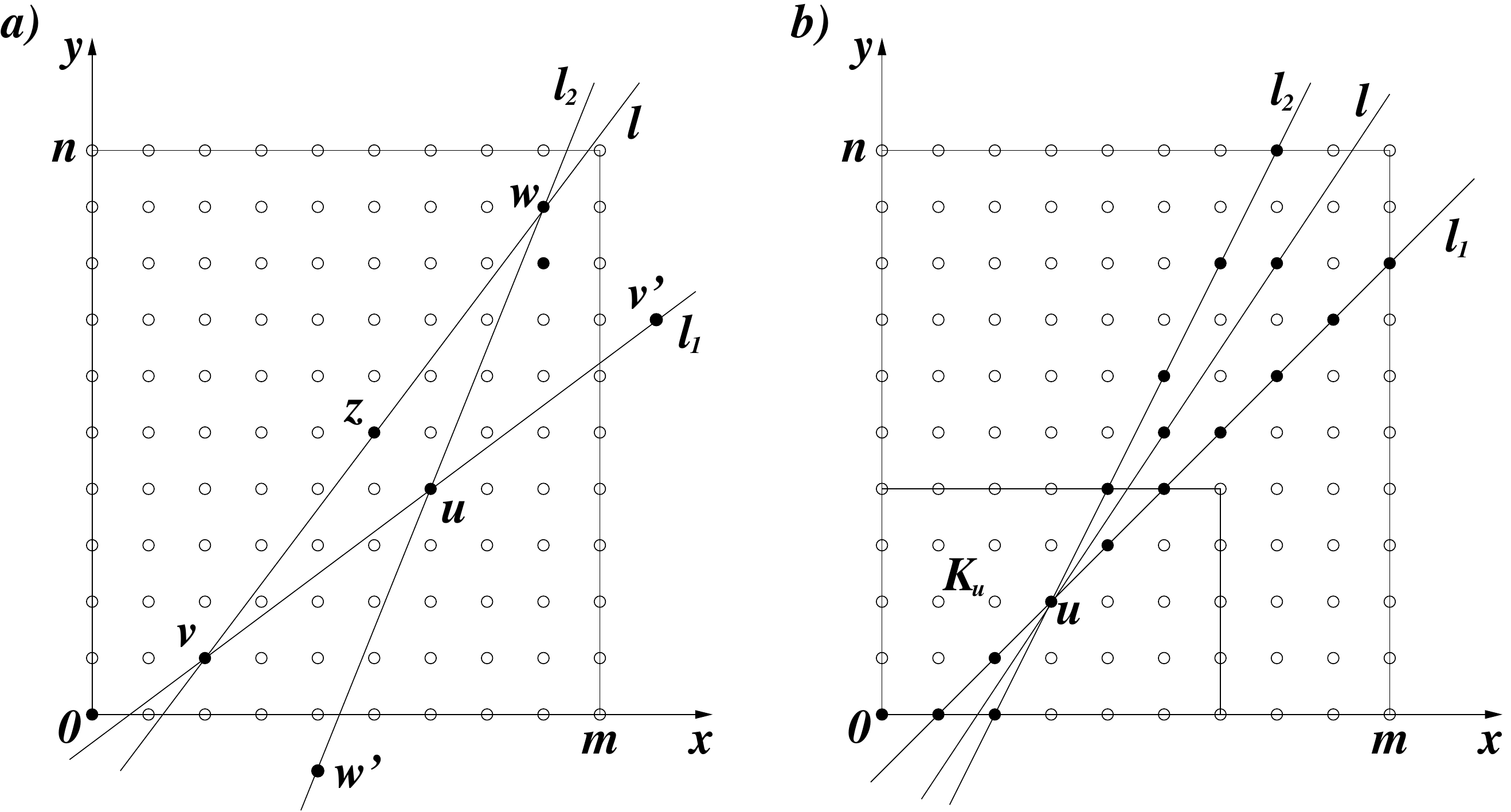}
\caption{a) The point $(w_x,w_y-1)$ belongs to $M(\l_1)\setminus M(\l_2)$, 
implying that $\l_1$ and $\l_2$ are not equivalent.
b) The line $\l$ lies between the adjacent $u$-stable lines $\l_1$ and $\l_2$.
The line $\l$ is stable but not $u$-stable since it has no integer points within the rectangle $K_u$, except $u$.
}
\label{Fig1}
\end{center}
\end{figure}

\begin{proof} Let $\l_1\cap \l_2=u\in\K_0$. Suppose that the line $\l_1$ passes through points $u\ne v\in\K_0$,
while the line $\l_2$ passes through points $u\ne w\in\K_0$.

If the line $\l_1$ is vertical, then $u_x=v_x=c$. 
In this case the line $\l_2$ cannot be horizontal, since otherwise the corresponding threshold function would be the zero function which is not in $\M$.
From $w\in M(\l_2)=M(\l_1)$ it follows that $w_x\leq c$. 
If $w_x<c$, then $(c,w_y)\in M(\l_1)\setminus M(\l_2)$ which contradicts $\l_1\sim \l_2$. 
Therefore, $w_x=c$ and $\l_2=\l_1$.

The other cases with a horizontal or vertical line are considered similarly. 
Assume that both lines $\l_1$, $l_2$ are inclined.

If the point $u$ lies on the border of $\K$, then
either $v\not\in M(\l_2)$ or $w\not\in M(\l_1)$, a contradiction to $\l_1\sim \l_2$.
Hence, $u$ is an internal point of the rectangle $\K$.

It is easy to see that if the lines $\l_1$ and $\l_2$ have opposite signs, 
then the line $x=u_x$ or $y=u_y$ contains a point from $M(\l_1)\bigtriangleup M(\l_2)$.

The remaining case to consider is the lines $\l_1$ and $\l_2$ having the same sign.
Without loss of generality, assume that they are positive. Then each of them crosses the left or bottom side of the rectangle $\K$.
Note that the case, when one line crosses the left side while the other crosses the bottom side, is
impossible since it would imply $(0,n)\in M(\l_1)\bigtriangleup M(\l_2)$ 
which contradicts $\l_1\sim \l_2$.
Without loss of generality, assume that both lines cross the bottom side of $\K$.

Suppose that the slopes of $\l_1$ and $\l_2$ are not equal. Without loss of generality, assume that
the slope of $\l_1$ is less than the slope of $\l_2$. Then $v_x<u_x<w_x$ and $v_y<u_y<w_y$. Hence,
line $\l$, passing through the points $v$ and $w$, is positive (Fig.~\ref{Fig1}a).

Let $z=\frac{v+w}{2}$ be the middle point of the interval $[v,w]$. Consider
the point $v'$ symmetric to $v$ with respect to the point $u$. It is easy to see that $v'\not\in\K_0$,
since otherwise $v'\in M(\l_1)\setminus M(\l_2)$ contradicts $\l_1\sim \l_2$.
In particular, we have $v'_x=2u_x-v_x>w_x$, which is equivalent to $u_x>z_x$.
Similarly, for the point $w'$ symmetric to $w$ with respect to $u$, we have $u_y<z_y$.
Since $u$ is an integral point and $z$ is the middle point of the interval with integral endpoints, 
a stronger inequality $z_y-u_y\geq\frac{1}{2}$ holds.

Consider an equation of the line $\l_1$ in the form
$y_1(x)=\frac{u_y-v_y}{u_x-v_x}(x-v_x)+v_y$ 
and an equation of the line $\l$ in the form
$y(x)=\frac{w_y-v_y}{w_x-v_x}(x-v_x)+v_y$. Define a function
$$
f(x)\eqdef y(x+v_x)-y_1(x+v_x)=\ls\frac{w_y-v_y}{w_x-v_x}-\frac{u_y-v_y}{u_x-v_x}\rs x.
$$
Trivially, we have $f(kx)=kf(x)$ for all $k$.

Since the line $\l_1$ is positive, and $z_x<u_x$, then $y_1(z_x)<y_1(u_x)=u_y$.
Hence, $f\ls\frac{w_x-v_x}{2}\rs=z_y-y_1(z_x)>z_y-u_y\geq\frac{1}{2}$.
Linearity of $f(x)$ implies
$$
w_y-y_1(w_x)=f(w_x-v_x)=2 f\ls\frac{w_x-v_x}{2}\rs > 1,
$$
which is equivalent to $y_1(w_x)<w_y-1$. Hence, $(w_x,w_y-1)\in M(\l_1)\setminus M(\l_2)$
which contradicts $\l_1\sim \l_2$.
The contradiction proves that the slopes of $\l_1$ and $\l_2$ are equal.
Since $\l_1$ and $\l_2$ have the common point $u$, Lemma~\ref{uK0} implies $\l_1=\l_2$.
\end{proof}

\begin{Lemma}\label{u2}
If lines $\l_1$ and $\l_2$ are stable and $\l_1\sim \l_2$, then $\l_1=\l_2$.
\end{Lemma}

\begin{proof} Assume that lines $\l_1$ and $\l_2$ are stable and $\l_1\sim \l_2$,
but $\l_1\ne \l_2$. By Lemma \ref{l2q}, $\l_1\cap \l_2=q\in\K$.

If $q\in\K_0$, then Lemma \ref{uK0} implies $\l_1=\l_2$, a contradiction.

Assume $q\not\in\K_0$. Let points $u,v\in\K_0$ be the closest to $q$
lying on the lines $\l_1$ and $\l_2$ respectively.
By Lemma~\ref{l12l0} there exists a stable line $\l_0$ equivalent to $\l_1$ such that $u\in \l_0$.
Lemma~\ref{uK0} applied to the lines $\l_0$ and $\l_1$ implies $\l_0=\l_1$,
which is impossible since $\l_0$ differs from $\l_1$ by construction.
This contradiction completes the proof.
\end{proof}

\begin{Lemma}\label{u1} Let the line $\l_1$ pass through a point $u\in\K_0$. If line $\l_2$
is stable and $\l_1\sim \l_2$, then $\l_2$ passes through $u$ as well.
\end{Lemma}

\begin{proof} If the line $\l_1$ is stable, then the statement immediately follows
from Lemma~\ref{uK0}. Hence, suppose that the line $\l_1$ is not stable, i.e., $\l_1\cap\K_0=\{u\}$.

According to Lemma~\ref{l2q}, $\l_1\cap \l_2=q\in\K$. Let $v\in\K_0$ be a point on the line $\l_2$ closest to the point $q$.

If $q\in\K_0$, i.e., $v=u=q$, then the lemma is proved.

Assume $q\not\in\K_0$. By Lemma \ref{l12l0}, there exists
a stable line $\l_0$ such that $u\in \l_0$ and $\l_0\sim \l_1$.
Applying Lemma \ref{u2} to the lines $\l_0$ and $\l_2$, we conclude that
$\l_0=\l_2$ and thus $u\in \l_2$.
\end{proof}

\begin{Def} A threshold function $f\in\M$ is called \emph{stable}, 
if there exists a stable line that defines $f$.
If there is no such line, the function $f$ is called \emph{unstable}.
\end{Def}

Examples of stable and unstable functions are given in Fig.~\ref{Fig0}.


\subsection{Number of unstable functions}

\begin{Def} Lemma~\ref{l12l0} implies that for an unstable function $f$, 
every line from $\L$ defining $f$ passes through the very same point from $K_0$.
We call this point the \emph{vertex} of $f$.
\end{Def}

\begin{Def} For $u\in \K_0$, let $\K_u$ be the largest axis-parallel rectangle 
contained in $\K$ with the center at $u$. Further let 
$\overline{\K_u}=\K\setminus\K_u$. 
Denote by $L_u$ the set of all lines passing through $u$ and by $S_u$ the set of 
all stable lines from $L_u$.
A line $\l$ is called \emph{$u$-stable} if $\l\in S_u$ and it passes through 
a point from $(\K_u\cap\K_0)\setminus \{u\}$.
\end{Def}

\begin{Def} Let $u\in\K_0$. Lines $\l_1,\ \l_2\in S_u$ are called
\emph{adjacent} if in the shortest rotation of a line about the point $u$ from the position $\l_1$
to the position $\l_2$, it does not meet any other lines from $S_u$.
Any intermediate line in this rotation is said to be \emph{lying between} the lines $\l_1$ and $\l_2$ (Fig.~\ref{Fig1}b).
\end{Def}

\begin{Lemma}\label{uss} If a line $\l$ passing through $u\in\K_0$ lies between
adjacent lines $\l_1,\ \l_2\in S_u$, and $\l\sim \l'$ for some stable line $\l'$,
then either $\l'=\l_1$ or $\l'=\l_2$.
\end{Lemma}
\begin{proof} Assume that $\l'$ differs from $\l_1$ and $\l_2$.
Let $v_1$ and $v_2$ be points from $\K_0$ different from $u$ lying
on the lines $\l_1$ and $\l_2$ respectively. Then Lemma~\ref{u1} implies that $u\in \l'$.
In the shortest rotation about the point $u$ towards the line $\l'$, the line $\l$
necessarily meets $\l_t$ for $t=1$ or $2$. But then the point $v_t$ belongs to
$M(\l)\bigtriangleup M(\l')$ and thus $\l$ and $\l'$ cannot be equivalent.
This contradiction completes the proof.
\end{proof}

For $u=\0$, there exists a unique unstable threshold 
function with the vertex at $u$. 
Namely, this function takes value $0$ only at the point $\0$
and thus is defined by any negative singular line.

If $u\ne\0$ is a corner vertex of the rectangle $\K$, 
then there is no unstable threshold function with the vertex at $u$. 
It is easy to see that any line passing through $u$ can be rotated about $u$ into 
an equivalent stable line.

Denote by $\K'_0$ the set $\K_0$ with excluded the corner vertices of the rectangle $\K$.

\begin{Lemma}\label{knuf} Let $u\in\K'_0$, and $\l_0$ be a line passing through
the points $\0$ and $u$. Then the number of unstable threshold functions with the
vertex at $u$ is equal to the number of $u$-stable lines, 
if the line $\l_0$ is $u$-stable; and one less otherwise.
\end{Lemma}
\begin{proof} Suppose that a line $\l\in L_u\setminus S_u$ lies between 
adjacent lines $\l_1,\ \l_2\in S_u$.
Note that $\l$ cannot be equivalent to any $u$-stable line, since the latter
contains a pair of symmetric (with respect to $u$) points, one of
which does not belong to $M(\l)$. Hence, if both lines $\l_1$ and $\l_2$ are
$u$-stable, then $\l$ is not equivalent to either of them
and thus by Lemma~\ref{uss} defines an unstable threshold function.

Therefore, the statement is true when $\K_u=\K$ (e.g., $u$ is the center of
the rectangle $\K$). For the rest of the proof, we assume that $u$ 
is not the center of $\K$.

Let us split the set $L_u$ into two:
\begin{eqnarray*}
L'_u & \eqdef & \{\l\in L_u\mid \l\cap\overline{\K_u}=\emptyset\};\\
L^{\prime\prime}_u & \eqdef & \{\l\in L_u\mid \l\cap\overline{\K_u}\ne\emptyset\}.
\end{eqnarray*}

Any unstable line from $L'_u$ lies between two adjacent stable lines from
$L'_u$, which are trivially $u$-stable as well.
Hence, the number of distinct unstable threshold functions defined 
by lines from $L'_u$ equals the number of $u$-stable lines in $L'_u$ minus 1.

Any unstable line $\l\in L^{\prime\prime}_u$ consists of two rays starting at $u$
such that one of them crosses $\overline{\K_u}$. 
Consider a rotation of $\l$ about $u$ such this ray moves towards the origin, following the shortest arc. 
Without loss of generality, we assume that in this rotation, the line $\l$ first meets the stable line $\l_1$. 
Then $\l$ cannot be equivalent to the line $\l_2$, since $\l_2$ is either $u$-stable or passes 
through a point from $\overline{\K_u}\cap\K_0$ not belonging to $M(\l)$.
On the other hand, $\l$ is equivalent to $\l_1$ if and only if $\l_1$ is not $u$-stable. 
Hence, the counting of unstable threshold functions with the vertex at $u$ corresponds to counting of $u$-stable lines.
The number of unstable threshold functions defined by lines from $L^{\prime\prime}_u$ equals 
the number of $u$-stable lines in $L^{\prime\prime}_u$.

Note that if the line $\l_0$ is $u$-stable, it is counted two times. 
In this case, the number of unstable threshold functions defined by lines from $L^{\prime\prime}_u$
is greater by $1$ as compared to the number of $u$-stable lines in $L^{\prime\prime}_u$.

Therefore, the total number of unstable threshold functions with the vertex at $u$ 
is equal to the number of $u$-stable lines if the line $\l_0$ is $u$-stable; and is less by $1$ otherwise.
\end{proof}

\begin{Def} Let
$$
U(p,q) \eqdef \card \{ (a,b)\in\Z^2\mid 1\leq a\leq p,\ 1\leq b\leq q,\ \gcd(a,b)=1\}.
$$
In other words,
\beq\label{Udef}
U(p,q) = \mathop{\sum_{i=1}^p\sum_{j=1}^q}\limits_{\gcd(i,j)=1} 1.
\eeq
\end{Def}

It is easy to see that $U(m,n)$ gives the number of singular inclined stable lines, each defined by 
the slope $\frac{b}{a}$ (i.e., passing through the points $\0$ and $(a,b)$) with $1\leq a\leq m$, $1\leq b\leq n$, and $\gcd(a,b)=1$.
We also remark that the values $U(k,k)$ are related to the probability of 
two random integers from $[1,k]$ being co-prime (see sequence \emph{A018805} in \cite{OEIS}).

\begin{Lemma}\label{l0nu} 
The number of points $u\in\K'_0$ such that the line  passing through the points $\0$ and $u$ is not $u$-stable equals $U(m,n)-1$.
\end{Lemma}
\begin{proof}
If a line $\l$ passes through exactly $k+1$ points $\0=u_0, u_1, \dots, u_k\in\K_0$ 
listed in the order of increasing distance from the point $\0$,
then $\l$ is $u_i$-stable for $i=1,2,\dots,k-1$ and is not $u_k$-stable.
Therefore, on each inclined singular stable line $\l$ there is exactly one point $u\in\K_0$ such that $\l$ is not $u$-stable.
Hence, the number of such points $u\in\K'_0$ equals the number of inclined singular stable lines, 
excluding the line $y=\frac{n}{m}x$ with $u=(m,n)\not\in\K'_0$, which is $U(m,n)-1$.
\end{proof}

\begin{Def} For real numbers $t$ and $k$, define a function
\beq\label{Vdef}
V(t,k) \eqdef 
\mathop{\sum_{i=1}^{\ceil{t}}\sum_{j=1}^{\ceil{k}}}\limits_{\gcd(i,j)=1} (t+1-i)(k+1-j).
\eeq
\end{Def}

\begin{Lemma}\label{Vprop}
For all $t\in\R, k\in\N$ the following equality holds
\begin{equation}\label{Vprop1}
V(t,k-1) + V(t,k) = 2 V\ls t, k-\frac{1}{2}\rs.
\end{equation}
For all $t, k\in\N$ the following equalities hold
\begin{equation}\label{Vprop2}
V(t,k) = \sum_{p=1}^t\sum_{q=1}^k U(p,q);
\end{equation}
\begin{equation}\label{Vprop3}
\sum_{p=1}^t U(p,k) = V(t,k) - V(t, k-1);
\end{equation}
\begin{equation}\label{Vprop4}
U(t,k) = V(t,k) - V(t,k-1) - V(t-1,k) + V(t-1,k-1);
\end{equation}
\begin{equation}\label{Vprop5}
U(t,k) + 2\ls V(t,k-1) + V(t-1,k)\rs = 4 V\ls t-\frac{1}{2},k-\frac{1}{2}\rs.
\end{equation}
\end{Lemma}
\begin{proof} Let $t\in\R, k\in\N$. To prove \eqref{Vprop1}, we use formula \eqref{Vdef}
\begin{eqnarray*}
V(t,k-1) + V(t,k) 
& = & \mathop{\sum_{i=1}^{\ceil{t}}\sum_{j=1}^{k-1}}\limits_{\gcd(i,j)=1}
(t+1-i)(k-j) 
+ \mathop{\sum_{i=1}^{\ceil{t}}\sum_{j=1}^k}\limits_{\gcd(i,j)=1} (t+1-i)(k+1-j) \\
& = & \mathop{\sum_{i=1}^{\ceil{t}}\sum_{j=1}^k}\limits_{\gcd(i,j)=1}(t+1-i)(2k+1-2j) 
= 2 \mathop{\sum_{i=1}^{\ceil{t}}\sum_{j=1}^k}\limits_{\gcd(i,j)=1}
(t+1-i)(k+\frac{1}{2}-j) \\
& = & 2V\ls t,k-\frac{1}{2}\rs.
\end{eqnarray*}

To prove \eqref{Vprop2}, we use formulae \eqref{Udef} and \eqref{Vdef}
\begin{eqnarray*}
\sum_{p=1}^t\sum_{q=1}^k U(p,q) 
& = & \sum_{p=1}^t\sum_{q=1}^k\mathop{\sum_{i=1}^p\sum_{j=1}^q}\limits_{\gcd(i,j)=1} 1
= \mathop{\sum_{i=1}^t\sum_{j=1}^k}\limits_{\gcd(i,j)=1} \sum_{p=i}^t
\sum_{q=j}^k 1 
= \mathop{\sum_{i=1}^t\sum_{j=1}^k}\limits_{\gcd(i,j)=1} (t+1-i)(k+1-j) \\
& = & V(t,k).
\end{eqnarray*}

To prove \eqref{Vprop3}, we use formula \eqref{Vprop2}
$$
V(t,k) - V(t, k-1) = \sum_{p=1}^t\sum_{q=1}^k U(p,q) -
\sum_{p=1}^t\sum_{q=1}^{k-1} U(p,q) = \sum_{p=1}^t U(p,k).
$$

Using \eqref{Vprop3}, we prove \eqref{Vprop4}
$$
U(t,k) = \sum_{p=1}^t U(p,k) - \sum_{p=1}^{t-1} U(p,k) = V(t,k) - V(t, k-1) -
V(t-1,k) + V(t-1, k-1).
$$

Finally, to prove equality \eqref{Vprop5}, we use formulae \eqref{Vprop4} and \eqref{Vprop1}
\begin{eqnarray*}
U(t,k) + 2\ls V(t,k-1) + V(t-1,k)\rs 
& = & V(t,k) + V(t,k-1) + V(t-1,k) + V(t-1,k-1) \\
& = & 2V\ls t,k-\frac{1}{2}\rs + 2V\ls t-1,k-\frac{1}{2}\rs \\
& = & 4 V\ls t-\frac{1}{2},k-\frac{1}{2}\rs.
\end{eqnarray*}
\end{proof}

\begin{Th}\label{knf} The number of unstable threshold functions in $\M$ is
$$
2mn - U(m,n) + 8 V\ls\frac{m-1}{2},\frac{n-1}{2}\rs.
$$
\end{Th}
\begin{proof}
We first notice that if $u_x\leq\frac{m}{2}$ and $u_y\leq\frac{n}{2}$, 
then the number of $u$-stable inclined lines equals $2U(u_x,u_y)$.
If $u$ lies on a side of the rectangle $\K$, then the line containing this side is the only $u$-stable line.
Otherwise, if $u$ is an internal point of $\K$, then both vertical and horizontal as well as
inclined lines passing through $u$ are $u$-stable.

Let us count the number of all $u$-stable lines for $u\in\K'_0$. 
Despite that counting depends on the parity of the integers $m$ and $n$, 
we will show that the result in all cases equals
\beq\label{uuK}
2(mn-1) + 8 V\ls\frac{m-1}{2},\frac{n-1}{2}\rs.
\eeq

Below we consider the different possible parities of $m$ and $n$.

If both $m$ and $n$ are odd, then \eqref{Vprop2} and Lemma~\ref{knuf} imply that the number of all $u$-stable lines for
$u=(i,j)\in [0,\frac{m-1}{2}]\times[0,\frac{n-1}{2}]\setminus \{\0\}$ equals
$$
\sum_{i=1}^{\frac{m-1}{2}} 1 + \sum_{j=1}^{\frac{n-1}{2}} 1 +
\sum_{i=1}^{\frac{m-1}{2}}\sum_{j=1}^{\frac{n-1}{2}} \ls 2U(i,j)+2\rs =
\frac{mn-1}{2} + 2 V\ls\frac{m-1}{2},\frac{n-1}{2}\rs.
$$
Due to the symmetry, the total number of $u$-stable lines for $u\in\K'_0$ is
four times as many.

Now let $m$ be even and $n$ be odd. Then the number of all $u$-stable lines for
$u=(i,j)\in [0,\frac{m}{2}-1]\times[0,\frac{n-1}{2}]\setminus \{\0\}$ equals
$$
\sum_{i=1}^{\frac{m}{2}-1} 1 + \sum_{j=1}^{\frac{n-1}{2}} 1 +
\sum_{i=1}^{\frac{m}{2}-1}\sum_{j=1}^{\frac{n-1}{2}} \ls 2U(i,j)+2\rs =
\frac{mn-n-1}{2} + 2 V\ls \frac{m}{2}-1, \frac{n-1}{2} \rs.
$$
Quadruplicated this number counts the $u$-stable lines for $u\in\K'_0$, 
except for $u$ lying on the line $x=\frac{m}{2}$. 
Due to property \eqref{Vprop3}, the $u$-stable lines for $u$ lying on the line $x=\frac{m}{2}$ 
can be counted as
$$
2\ls 1 + \sum_{j=1}^{\frac{n-1}{2}} \ls 2U(\frac{m}{2},j)+2\rs\rs = 2n + 4\ls
V\ls\frac{m}{2},\frac{n-1}{2}\rs - V\ls\frac{m}{2}-1,\frac{n-1}{2}\rs\rs.
$$
Hence, the total number of $u$-stable lines for $u\in\K'_0$ equals
$$
\bea{l}
4\ls \frac{mn-n-1}{2} + 2 V\ls \frac{m}{2}-1, \frac{n-1}{2} \rs \rs +
2n + 4\ls V\ls\frac{m}{2},\frac{n-1}{2}\rs - V\ls\frac{m}{2}-1,\frac{n-1}{2}\rs\rs
\\
= 2(mn-1) + 4 \ls V\ls\frac{m}{2},\frac{n-1}{2}\rs +
V\ls\frac{m}{2}-1,\frac{n-1}{2}\rs\rs = 2(mn-1) + 8V\ls\frac{m-1}{2},\frac{n-1}{2}\rs.
\eea
$$
Here we used the property \eqref{Vprop1}.

The case of odd $m$ and even $n$ is considered similarly.

Finally, let both $m$ and $n$ be even. The number of all $u$-stable lines for
$u=(i,j)\in [0,\frac{m}{2}-1]\times[0,\frac{n}{2}-1]\setminus \{\0\}$ equals
$$
\sum_{i=1}^{\frac{m}{2}-1} 1 + \sum_{j=1}^{\frac{n}{2}-1} 1 +
\sum_{i=1}^{\frac{m}{2}-1}\sum_{j=1}^{\frac{n}{2}-1} \ls 2U(i,j)+2\rs =
\frac{mn-m-n}{2} + 2 V\ls \frac{m}{2}-1, \frac{n}{2}-1 \rs.
$$
Quadruplicated this number equals the number of all $u$-stable lines for $u\in\K'_0$, 
except for those lying on the lines $x=\frac{m}{2}$ or $y=\frac{n}{2}$.
Properties \eqref{Vprop3} and \eqref{Vprop5} imply that
the number of $u$-stable lines for $u$ lying on these lines equals
$$
\bea{l}
2\ls 1 + 1 + \sum_{i=1}^{\frac{m}{2}-1} \ls 2U(i,\frac{n}{2})+2\rs +
\sum_{j=1}^{\frac{n}{2}-1} \ls 2U(\frac{m}{2},j)+2\rs\rs +
2U\ls\frac{m}{2},\frac{n}{2}\rs+2 \\
= 2(m+n-1) + 2U\ls\frac{m}{2},\frac{n}{2}\rs + 
4\ls V\ls\frac{m}{2}-1,\frac{n}{2}\rs + V\ls\frac{m}{2},\frac{n}{2}-1\rs -
2V\ls\frac{m}{2}-1,\frac{n}{2}-1\rs\rs \\
= 2(m+n-1) + 2V\ls\frac{m-1}{2},\frac{n-1}{2}\rs - 8V\ls\frac{m}{2}-1,\frac{n}{2}-1\rs.
\eea
$$
Summing up the results, we get \eqref{uuK}.

Lemmas \ref{knuf} and \ref{l0nu} imply that the total number of unstable threshold
functions with the vertex in $\K'_0$ is equal to \eqref{uuK} minus $U(m,n)-1$.
Finally, we need to add $1$ for the single unstable threshold function with the vertex at $\0$.
\end{proof}


\subsection{Number of stable threshold functions}

\begin{Th}\label{kuf} The number of stable threshold functions in $\M$ equals
$$
m + n + U(m,n) + 2 V(m,n) - 8V\ls\frac{m-1}{2},\frac{n-1}{2}\rs.
$$
\end{Th}
\begin{proof}
Consider any stable line $\l$ of a positive slope passing through $\0$.
Let a point $(a,b)\in l\cap\K_0$ be the closest to $\0$, implying that $\gcd(a,b)=1$.

Consider all stable lines parallel to $\l$. Every such line is defined 
by a pair of points $(x,y), (x+a,y+b)\in\K_0$ on it, where $(x,y)\in\K_0$ 
is the closest point to $\0$. 
Such pairs are uniquely defined by the following constraints
$$
\begin{cases}
x<a\ \text{or}\ y<b,\\
x+a\leq m,\\
y+b\leq n.
\end{cases}
$$
Let $P$ be the set of all points $(x,y)\in\K_0$ satisfying these constraints.
Then the stable lines parallel to $\l$ and the elements of $P$ are in 
one-to-one correspondence.

If $a>\frac{m}{2}$ or $b>\frac{n}{2}$, then the set $P$ equals
$P_1\eqdef\{(x,y)\in\K_0\mid 0\leq x\leq m-a,\ 0\leq y\leq n-b\}$. In this
case the number of stable lines parallel to $\l$ is $(m+1-a)(n+1-b)$.

If $a\leq\frac{m}{2}$ and $b\leq\frac{n}{2}$, then the set $P$ equals
$$
P_1 \setminus
\{(x,y)\in\K_0\mid a\leq x\leq m-a,\ b\leq y\leq n-b\}.
$$
Hence, in this case the number of stable lines parallel to $\l$ is less 
than before by
$$
(m+1-2a)(n+1-2b) = 4\ls\frac{m-1}{2}+1-a\rs\ls\frac{n-1}{2}+1-b\rs .
$$

Summing over all pairs $(a,b)$ fulfilling the constraints and using formula 
\eqref{Vdef}, we derive that the total number of stable 
lines of a positive (negative) slope is
$V(m,n)-4V\ls\frac{m-1}{2},\frac{n-1}{2}\rs$,
while the total number of inclined stable lines is twice as many.

Since each inclined stable line passing through $\0$ defines two distinct threshold
functions, we further increase the count by the number of such lines, i.e., $U(m,n)$.

Finally, we take into account $m$ vertical lines $x=i$ for $i=\overline{0,m-1}$ 
and $n$ horizontal lines $y=j$ for $j=\overline{0,n-1}$ to complete the proof.
\end{proof}

Adding up the results of Theorems~\ref{knf} and \ref{kuf} and noticing that $N(m,n)=2(|\M|+1)$ completes
the proof of Theorem~\ref{Total}.


\section{Asymptotics of $N(m,n)$}\label{Sasympt}

In this section we prove the following theorem, which together with \eqref{FormulaPN} implies asymptotics \eqref{AsympPk2}.

\begin{Th}\label{asympN} For $m\geq n$, the following asymptotics holds
$$
N(m,n) = \frac{6}{\pi^2} m^2n^2 + O(m^2n\ln n)
$$
$$
N(m,n) = 2((n+1)\Psi(n) - \Phi(n))m^2 + O(mn^3)
$$
where $\varphi(t)$ is the totient function and
\begin{itemize}
\item $\Phi(k)\eqdef\sum_{i=1}^k \varphi(i) = \frac{3}{\pi^2} k^2 +
O(k\ln k)$ (Dirichlet's Theorem, see \cite{Prachar});
\item $\Psi(k)\eqdef\sum_{i=1}^k \frac{\varphi(i)}{i} = \frac{6}{\pi^2} k + O(\ln k)$ (see \cite{Prachar}).
\end{itemize}
\end{Th}

We remark that the first of the two asymptotics for $N(m,n)$ given in 
Theorem~\ref{asympN} is more suitable for $m$ and $n$ of the same magnitude, 
while the second asymptotics is better when $m\gg n$.

\begin{Lemma}\label{OOO} Let $k$ be a positive integer and $s\geq 0$. Then
$$
\sum_{t=1}^k\frac{1}{t}=O(\ln k);
$$
$$
\sum_{t=k+1}^{\infty}\frac{1}{t^s} = O\ls\frac{1}{k^{s-1}}\rs,\quad (s\ne 1);
$$
$$
\sum_{t=1}^k t^s\ln t = O\ls k^{s+1}\ln k\rs;
$$
$$
\sum_{t=1}^k t^s = \frac{k^{s+1}}{s+1} + O(k^s) = O(k^{s+1}).
$$
\end{Lemma}

\begin{proof} The statement follows from integral estimates of the sums.
\end{proof}

\begin{Th}\label{umk} For $t\geq k$, the following inequality holds
$$
U(t,k) = U(k,t) = \frac{6}{\pi^2}tk + O(t\ln k);
$$
\end{Th}
\begin{proof} Note that there exist exactly $\floor{m/p}$ positive integers
not exceeding $m$ that are divisible by $p$. Hence, there are
$\floor{t/p}\floor{k/p}$ pairs $(a,b)$, $1\leq a\leq t$, $1\leq b\leq k$, 
whose greatest common divisor is divisible by $p$.
The inclusion-exclusion principle \cite{Stanley} for the number of pairs $(a,b)$
with $\gcd(a,b)$ not divisible by any prime $p$ (i.e., $\gcd(a,b)=1$) 
gives an exact formula
$$
U(t,k) = \sum_{s=1}^{\infty} \mu(s) \floor{\frac{t}{s}}\floor{\frac{k}{s}},
$$
where $\mu(s)$ is the M\"obeus function.

We approximate $U(t,k)$ with the function
$$
\hat U(t,k)\eqdef \sum_{s=1}^{\infty} \mu(s) \frac{t}{s}\frac{k}{s} =
tk\sum_{s=1}^{\infty} \frac{\mu(s)}{s^2} = tk\frac{6}{\pi^2}.
$$
and bound the absolute value of the difference $U(t,k)-\hat U(t,k)$ 
as follows
$$
\left|U(t,k)-\hat U(t,k)\right| = \left|\sum_{s=1}^{\infty} \mu(s)
\ls\floor{\frac{t}{s}}\floor{\frac{k}{s}} - \frac{t}{s}\frac{k}{s}\rs\right| \leq
\sum_{s=1}^{\infty} d(t,k,s),
$$
where
$$
d(t,k,s)\eqdef \frac{tk}{s^2} - \floor{\frac{t}{s}}\floor{\frac{k}{s}}.
$$

For $s>k$, we have $\floor{\frac{k}{s}}=0$ and thus $d(t,k,s)=\frac{tk}{s^2}$. 
For $s\leq k$, we bound $d(t,k,s)$ as follows:
$$
d(t,k,s) <
\frac{tk}{s^2} - \ls\frac{t}{s}-1\rs\ls\frac{k}{s}-1\rs
= \frac{t+k}{s} - 1 < \frac{2t}{s}.
$$

Applying Lemma~\ref{OOO}, we finally have
$$
\left|U(t,k)-\hat U(t,k)\right| < \sum_{s=1}^{\infty} d(t,k,s) <
\sum_{s=1}^k \frac{2t}{s} + \sum_{s=k+1}^{\infty} \frac{tk}{s^2} =
O(t\ln k).
$$
\end{proof}

We remark that Dirichlet's Theorem is a particular case of Theorem \ref{umk} for $t=k$.

\begin{Th}\label{vmk} For $t\geq k$, the following asymptotics holds
$$
V(t,k) = V(k,t) = \frac{3}{2\pi^2} t^2 k^2 + O(t^2 k\ln k).
$$
\end{Th}
\begin{proof} We use formula \eqref{Vprop2} and property $U(i,j)=U(j,i)$ as follows:
$$
V(t,k) = \sum_{i=1}^t\sum_{j=1}^k U(i,j) = 2\sum_{i=1}^k\sum_{j=1}^{i-1} U(i,j) +
\sum_{i=1}^k U(i,i) + \sum_{i=k+1}^t\sum_{j=1}^k U(i,j).
$$
According to Lemmas \ref{OOO} and \ref{umk},
$$
\sum_{i=1}^k\sum_{j=1}^{i-1} U(i,j) = \sum_{i=1}^k\sum_{j=1}^{i-1}
\ls\frac{6}{\pi^2}ij +
O(i\ln j)\rs = \frac{6}{\pi^2} \sum_{i=1}^k \ls i\sum_{j=1}^{i-1} j +
O\ls i \sum_{j=1}^{i-1}\ln j\rs\rs =
$$
$$
= \frac{6}{\pi^2}\sum_{i=1}^k \ls i\frac{i^2}{2} + O(i^2\ln i)\rs =
\frac{3}{\pi^2} \sum_{i=1}^k i^3 + O\ls\sum_{i=1}^k i^2\ln i\rs =
\frac{3}{4\pi^2} k^4 + O(k^3\ln k).
$$
Similarly,
$$
\sum_{i=1}^k U(i,i) = \sum_{i=1}^k O(i^2) = O(k^3);
$$
and
\begin{eqnarray*}
\sum_{i=k+1}^t\sum_{j=1}^k U(i,j) 
& = & \sum_{i=k+1}^t\sum_{j=1}^k \ls\frac{6}{\pi^2}ij + O(i\ln j)\rs 
= \sum_{i=k+1}^t \ls \frac{6}{\pi^2}i\frac{k^2}{2} + O(ik\ln k)\rs \\
& = & \frac{3}{2\pi^2}k^2(t^2-k^2) + O(t^2 k\ln k).
\end{eqnarray*}
Therefore,
\begin{eqnarray*}
V(t,k) & = & 2\ls \frac{3}{4\pi^2} k^4 + O(k^3\ln k)\rs + O(k^3) + \frac{3}{2\pi^2}k^2(t^2-k^2) + O(t^2 k\ln k) \\
& = & \frac{3}{2\pi^2} t^2 k^2 + O(t^2 k\ln k).
\end{eqnarray*}
\end{proof}

\begin{Th}\label{umkC} For $t\geq k$, the following asymptotics holds
$$
U(t,k) = U(k,t) = \Psi(k) t + O(k^2).
$$
\end{Th}
\begin{proof} We use formula \eqref{Udef}
\begin{eqnarray*}
U(k,t) 
& = & \mathop{\sum_{i=1}^k\sum_{j=1}^t}\limits_{\gcd(i,j)=1} 1 =
\sum_{i=1}^k\ls\sum_{s=0}^{\floor{\frac{t}{i}}-1}\sum_{j=si+1\atop\gcd(j,i)=1}^{si+i} 1 +
\sum_{j=\floor{t/i}i+1\atop\gcd(j,i)=1}^t 1\rs 
= \sum_{i=1}^k \floor{\frac{t}{i}}\varphi(i) + O(k^2) 
= \sum_{i=1}^k\frac{t}{i}\varphi(i) + O(k^2) \\
& = & \Psi(k) t + O(k^2).
\end{eqnarray*}
\end{proof}

\begin{Th}\label{vmkC} For $t\geq k$, the following asymptotics holds
$$
V(t,k) = V(k,t) = \frac{(k+1)\Psi(k) - \Phi(k)}{2}t^2 + O(tk^3).
$$
\end{Th}
\begin{proof} We use formula \eqref{Vdef}
\begin{eqnarray*}
V(k,t) 
& = & \mathop{\sum_{i=1}^k\sum_{j=1}^t}\limits_{\gcd(i,j)=1} (k+1-i)(t+1-j) \\
& = & \sum_{i=1}^k (k+1-i)\ls\sum_{s=0}^{\floor{\frac{t}{i}}-1}\sum_{j=1\atop\gcd(j,i)=1}^{i}
(t+1-si-j) + \sum_{j=1\atop\gcd(j,i)=1}^{t\bmod i} (t\bmod i + 1 - j)\rs.
\end{eqnarray*}
Neglecting terms of order $O(tk^3)$, we have
\begin{eqnarray*}
V(k,t) & = &
\sum_{i=1}^k (k+1-i)\sum_{s=0}^{\floor{\frac{t}{i}}-1}\sum_{j=1\atop\gcd(j,i)=1}^i
(t-si) + O(tk^3) = \frac{t^2}{2}\sum_{i=1}^k (k+1-i)\frac{\varphi(i)}{i} + O(tk^3) \\
& = & \frac{(k+1)\Psi(k) - \Phi(k)}{2}t^2 + O(tk^3).
\end{eqnarray*}
\end{proof}

Theorem \ref{asympN} now follows from Theorems \ref{Total}, \ref{vmk}, and \ref{vmkC}.

\section{Relation to teaching sets}\label{SecTeaching}

A \emph{teaching set} \cite{Goldman95,ShevZol98,ZolShev99} of a threshold function $f:\K_0\to\{0,1\}$ is
a subset $T\subset\K_0$ such that for any other threshold function $g\ne f$, there exists $t\in T$
with $g(t)\ne f(t)$. A teaching set of minimal cardinality is called \emph{minimum teaching set}.
A minimum teaching set of a two-dimensional threshold function consists of either 3 or 4 points \cite{ShevZol98}.


It is easy to see that any unstable threshold function $f\in\M$ has a teaching set of size 3. 
Indeed, let $u\in\K_0$ be the vertex of $f$ and $\l\ni u$ be a line defining $f$.
Then $\l$ lies between two adjacent stable lines $\l_1$ and $\l_2$ passing respectively through 
some points $u_1\ne u$ and $u_2\ne u$ with $f(u_1)=f(u_2)\ne f(u)$. 
Then $\{ u, u_1, u_2 \}$ forms a teaching set of $f$.

However, the size of a minimum teaching set $T$ of a stable threshold function $f\in\M$ may be $3$ or $4$.
Namely, $|T|=3$ if the complement threshold function $\tilde f(x,y)=1-f(m-x,n-y)$ is unstable;
and $|T|=4$ if $\tilde f(x,y)$ is stable. Therefore, the stable threshold functions can be partitioned 
into two classes, depending on the size of their minimal teaching sets. 
Unfortunately, our results do not allow to compute the number of threshold
functions in each class, which we pose as an open problem.


\section*{Acknowledgements}

The author thanks Nikolai Zolotykh and Valeriy Shevchenko for posing
the problem of counting threshold functions and related invaluable discussions, 
and Keith Conrad for reviewing an earlier version of the paper.





\bibliographystyle{siam}
\bibliography{porog_en.bib}

\end{document}